\documentclass[preprint,12pt]{elsarticle}
\usepackage{amssymb}
\usepackage{amsthm}
\usepackage{url}
\newdefinition{theorem}{Theorem}
\newdefinition{lemma}{Lemma}
\newdefinition{corollary}{Corollary}
\newdefinition{proposition}{Proposition}
\newdefinition{remark}{Remark}

\newcommand{\re}{\mathrm{Re}}

\begin{document}

\begin{frontmatter}

\title{On Some Explicit Formulas for Bernoulli Numbers and Polynomials}

\author{Lazhar Fekih-Ahmed}

\address{\'{E}cole Nationale d'Ing\'{e}nieurs de Tunis, BP 37, Le Belv\'{e}d\`{e}re 1002 , Tunis, Tunisia}
\ead{lazhar.fekihahmed@enit.rnu.tn}

\begin{abstract}
We provide direct elementary proofs of several explicit
expressions for Bernoulli numbers and Bernoulli polynomials. As a
byproduct of our method of proof, we provide natural definitions
for generalized Bernoulli numbers and polynomials of complex
order.

\end{abstract}


\end{frontmatter}


\section{Introduction}\label{sec1}

As Gould pointed out in his survey article \cite{gould}, a
well-known explicit formula  for Bernoulli numbers, which dates
back to Worpitzky \cite{worpitzky}, is given by the double sum

\begin{equation}\label{sec1-eq1}
B_{k}=\sum_{n=0}^{k}\frac{1}{n+1}\sum_{j=0}^{n}(-1)^{j}{n\choose
j}j^{k}, k\ge 0,
\end{equation}
where ${n\choose j}$ is the binomial coefficient. The first few
values of formula (\ref{sec1-eq1}) are $B_0=1$ $B_1=-1/2$,
$B_2=1/6$, $B_3=0$, $B_4=-1/30$, $B_5=0$, $B_6=1/42$, $B_7=0$,
$B_8=-1/30$ etc.

For the sake of convenience and to agree with our notation,  the
lower limits of summation in both sums in (\ref{sec1-eq1}) will be
changed so that the above sum is given by the following equivalent
form

\begin{equation}\label{sec1-eq2}
B_{k}=\sum_{n=1}^{k}\frac{1}{n+1}\sum_{j=1}^{n}(-1)^{j}{n\choose
j}j^{k}, k\ge 1.
\end{equation}

If we define the forward differences $\Delta_{n}(k)$ by

\begin{equation}\label{sec1-eq3}
\Delta_{n}(k)= \sum_{j=1}^{n}(-1)^{j}{n\choose j}j^{k},
\end{equation}

equation (\ref{sec1-eq2}) can be rewritten as

\begin{equation}\label{sec1-eq4}
B_{k}=\sum_{n=1}^{k}\frac{1}{n+1}\Delta_{n}(k), k\ge 1.
\end{equation}

In this note, our main result is  proving the following two
explicit formulas for Bernoulli numbers

\begin{equation}\label{sec1-eq5}
B_{k}=(-1)^{k+1}\sum_{n=1}^{k}\frac{1}{n(n+1)}\Delta_{n}(k), k\ge
1,
\end{equation}

and
\begin{equation}\label{sec1-eq6}
B_{k}=(-1)^{k+1}\sum_{n=1}^{k+1}\frac{1}{n^2}\Delta_{n}(k+1), k\ge
0.
\end{equation}

While formula (\ref{sec1-eq1}) is a well-known explicit formula
and several proofs have been provided by different authors, the
two formulas (\ref{sec1-eq5}) and (\ref{sec1-eq6}) are
less-known\footnote{Formula (\ref{sec1-eq5}) has been given in
\cite[formula (37)]{worpitzky}. The same formula is mentioned in
\cite[formula LXV on page 83]{Saalschutz}. Formula
(\ref{sec1-eq6}) is also mentioned in \cite[formula LXIII on page
82]{Saalschutz}. The proofs in \cite{Saalschutz} use the identity
$\Delta_{n}(k)=n(\Delta_{n}(k-1)-\Delta_{n-1}(k-1))$ and the
property that odd-indexed Bernoulli numbers have zero values. In
\cite[page 83]{Saalschutz}  the identity is written in terms of
$a_n^k$ as $a_{n}^k=n(a_{n}^{k-1}+a_{n-1}^{k-1})$ so that
$\Delta_{n}(k)=(-1)^{n}a_n^k=(-1)^{n}n!S(k,n)$, where $S(k,n)$ are
Stirling numbers of the second kind.}.

As a byproduct of our method of proof we provide extensions to
Bernoulli polynomials. We further provide natural definitions for
generalized Bernoulli numbers and polynomials of complex order.

\section{Fractional Derivatives}\label{sec2}

Suppose that the function $\psi(t)$ is holomorphic and that
$\lim_{t\to -\infty}\psi(t)=0$. According to Laurent
\cite{laurent}, the fractional derivative of order $\alpha \in
\mathbb{C}$ between the points $-\infty$ and $x \in \mathbb{R}$ of
the function $\psi(t)$ is given by the contour
integral\footnote{The derivative can of course be defined for all
$x \in \mathbb{C}$.}

\begin{equation}\label{sec2-eq1}
I^{\alpha}\psi(x)=\frac{\Gamma(1+\alpha)}{2\pi
i}\int_{\mathcal{C}}\frac{\psi(t)}{(t-x)^{\alpha +1}}\,dt,
\end{equation}

where $\mathcal{C}$ is the Hankel contour consisting of the three
parts $C=C_{-}\cup C_{\epsilon}\cup C_{+}$:  a path  which extends
from $(-\infty,-\epsilon)$, around the point $x$ counter clockwise
on a circle of center the point $x$ and of radius $\epsilon$ and
back to $(-\epsilon,-\infty)$, where $\epsilon$ is an arbitrarily
small positive number.

When $\re(1+\alpha)\ge 0$, there is no ambiguity in the definition
of $I^{\alpha}\psi(x)$. The integrals along $C_{-}$ and $C_{+}$
cancel each other. $I^{\alpha}\psi(x)$ is thus equal to the
integral along $C_{\epsilon}$ and this integral can be easily
evaluated by residue calculus. In particular, when $\alpha=0$,
formula (\ref{sec2-eq1}) is simply Cauchy's formula:

\begin{equation}\label{sec2-eq2}
I^{0}\psi(x)=\psi(x),
\end{equation}

and when $\alpha=n$ is a positive integer, Laurent's contour
integral $I^{\alpha}\psi(x)$  gives the classical derivative of
$\psi(t)$ at the point $x$:

\begin{equation}\label{sec2-eq3}
I^{n}\psi(x)= \psi^{(n)}(x).
\end{equation}

When $\re(1+\alpha)<1$, the portion of the integral along the
circle $C_{\epsilon}$ is zero. The integral along the remaining
portions of the contour is estimated using a proper determination
of the multi-valued function $(t-x)^{-\alpha-1}$. If we choose the
cut long the semi-axis $(-\infty,x)$, then
$(t-x)^{-\alpha-1}=e^{-(\alpha+1)(\log(t-x)-i\pi)}$ along $C_{-}$
and $(t-x)^{-\alpha-1}=e^{-(\alpha+1)(\log(t-x)+i\pi)}$ along
$C_{+}$, where $\log(t-x)$ is purely real when $t-x>0$. Moreover,
$t=re^{-i\pi}$ along $C_{-}$ and $t=re^{i\pi}$ along $C_{+}$, as
$r$ varies from $\epsilon$ to $+\infty$. The integral
(\ref{sec2-eq1}) becomes

\begin{equation}\label{sec2-eq4}
I^{\alpha}\psi(x)=\frac{\Gamma(1+\alpha)}{2\pi i} (e^{-\alpha \pi
i}-e^{\alpha \pi i})\int_{x}^{\infty}\psi(-r)
(r-x)^{-\alpha-1}\,dr.
\end{equation}

Finally, using the reflection formula of the Gamma function, we
obtain

\begin{equation}\label{sec2-eq5}
I^{\alpha}\psi(x)
=\frac{1}{\Gamma(-\alpha)}\int_{x}^{\infty}\psi(-r)
(r-x)^{-\alpha-1}\,dr.
\end{equation}

\section{Specializing to $\zeta(s)$}\label{sec3}

For the particular case $x=0$ and $\alpha =-s$, the fractional
derivative of order $-s$ at zero is given by

\begin{equation}\label{sec3-eq1}
I^{-s}\psi(0)=\frac{\Gamma(1-s)}{2\pi i}\int_{\mathcal{C}}\psi(t)
t^{s -1}\,dt
\end{equation}

when $\re(1-s)\ge 0$, and by

\begin{equation}\label{sec3-eq2}
I^{-s}\psi(0) =\frac{1}{\Gamma(s)}\int_{0}^{\infty}\psi(-t)
t^{s-1}\,dt
\end{equation}

when $\re(s)>0$.

Now by an appropriate choice of  $\psi(t)$, we will be able to
write the Riemann zeta function as the fractional derivative of
$\psi(t)$ at $0$. Indeed, let

\begin{equation}\label{sec3-eq3}
\psi(-t)=\phi(t)=\frac{d}{dt}\bigg(\frac{-t}{e^{t}-1}\bigg)=\frac{te^{-t}}{(1-e^{-t})^2}-\frac{e^{-t}}{1-e^{-t}}.
\end{equation}

We have shown in \cite{lazhar} that the Riemann zeta function has
the integral representation\footnote{The definition of $\phi(t)$
used here is different from the one in the cited paper.}

\begin{equation}\label{sec3-eq4}
(s-1)\zeta(s)
=\frac{1}{\Gamma(s)}\int_{0}^{\infty}\phi(t)t^{s-1}\,dt, \re(s)>0,
\end{equation}

and that for all $s\in \mathbb{C}$

\begin{equation}\label{sec3-eq5}
(s-1)\zeta(s)=\frac{\Gamma(1-s)}{2\pi
i}\int_{\mathcal{C}}\psi(t)t^{s-1}\,dt,
\end{equation}

where $\mathcal{C}$ is the same Hankel contour used in equation
(\ref{sec2-eq1}). Comparing with equations (\ref{sec3-eq1}) and
(\ref{sec2-eq2}), we easily obtain

\begin{equation}\label{sec3-eq6}
I^{-s}\big[\psi(t)\big]_{t=0}=(s-1)\zeta(s).
\end{equation}

That is, $(s-1)\zeta(s)$ is simply the fractional derivative of
order $-s$ of the function $\psi(t)$ at the origin. Furthermore,
for an integer $k\ge 2$, the derivative of order $k-1$ of
$\psi(t)$ at the point $t=0$ is, by Laurent's definition, given by

\begin{equation}\label{sec3-eq7}
I^{k-1}\big[\psi(t)\big]_{t=0}=\psi^{(k-1)}(0)=\frac{\Gamma(k)}{2\pi
i}\int_{\mathcal{C}}\psi(t) t^{-k}\,dt.
\end{equation}

Having achieved this, we know also that the Bernoulli numbers are
usually defined using the generating function

\begin{equation}\label{sec3-eq8}
\frac{t}{e^{t}-1}=\sum_{k=0}^{\infty}\frac{B_k}{k!}t^k,\quad
|t|<2\pi.
\end{equation}

These numbers can also be defined in terms of the function
$\psi(t)$ instead, since we have

\begin{equation}\label{sec3-eq9}
\psi(-t)=\frac{d}{dt}\bigg(\frac{-t}{e^{t}-1}\bigg)=\sum_{k=1}^{\infty}\frac{-B_k}{(k-1)!}t^{k-1}.
\end{equation}

Therefore,

\begin{eqnarray}
B_k&=&(-1)^{k}\psi^{(k-1)}(0)=(-1)^{k}I^{k-1}\big[\psi(t)\big]_{t=0}=(-1)^{k+1}k\zeta(1-k),\quad\rm{or}
\nonumber\\
B_k&=&-\phi^{(k-1)}(0)=-I^{k-1}\big[\phi(-t)\big]_{t=0}=k\zeta(1-k).\label{sec3-eq10}
\end{eqnarray}

This last equation is the basis of all our subsequent derivations.
It relates  the Bernoulli numbers, the Riemann-zeta function and
fractional derivatives in a single equation.

\section{The First Explicit Formula for $B_k$}\label{sec4}

In \cite{lazhar}, we have also obtained  a globally convergent
series representation of the Riemann zeta function. It is given by
the formula

\begin{equation}\label{sec4-eq1}
(s-1)\zeta(s)=\sum_{n=1}^{\infty}\frac{S_{n}(s)}{n+1},\qquad\rm{with}
\end{equation}

\begin{equation}\label{sec4-eq2}
S_{n}(s) = \sum_{k=0}^{n-1}(-1)^{k}{n-1\choose k}(k+1)^{-s}~{\rm
for}~ n\ge 2,
\end{equation}

and $S_{1}(s)=1$. We have also shown that when $\re(s)>0$ the sum
$S_{n}(s)$ can be rewritten as
\begin{equation}\label{sec4-eq3}
S_{n}(s)=\frac{1}{\Gamma(s)}\int_{0}^{\infty}(1-e^{-t})^{n-1}e^{-t}t^{s-1}\,dt,
\end{equation}

and that the function $\phi(t)$ verifies
\begin{equation}\label{sec4-eq4}
\phi(t)=\sum_{n=1}^{\infty}\frac{(1-e^{-t})^{n-1}e^{-t}}{n+1}
\end{equation}

uniformly\footnote{When $\re(s)<1$, $S_{n}(s)$ can obviously be
written as
\begin{equation}\nonumber
S_{n}(s)=\frac{\Gamma(1-s)}{2\pi
i}\int_{\mathcal{C}}(1-e^{t})^{n-1}e^{t}t^{s-1}\,dt.
\end{equation}
}  for $0<t<\infty$ .

Because of equation (\ref{sec4-eq3}), the definition of fractional
derivative (\ref{sec3-eq4}) and uniform convergence, we may
interchange summation and integration inside the integral
sign\footnote{This have been rigourously proved in
\cite{lazhar}.}. Thus, we may apply the fractional derivative
operator termwise to obtain

\begin{equation}\label{sec4-eq5}
I^{-s}\big[\psi(t)\big]_{t=0}=I^{-s}\big[\phi(-t)\big]_{t=0}=\sum_{n=1}^{\infty}\frac{I^{-s}\big[(1-e^{t})^{n-1}e^{t}\big]_{t=0}}{n+1}.
\end{equation}

Particularly, for $-s=k-1$ we obtain

\begin{eqnarray}
I^{k-1}\big[\phi(-t)\big]_{t=0}&=&\sum_{n=1}^{\infty}\frac{I^{k-1}\big[(1-e^{t})^{n-1}e^{t}\big]_{t=0}}{n+1}\nonumber\\
&=&\sum_{n=1}^{\infty}\frac{1}{n+1}\frac{d^{k-1}}{dt^{k-1}}\big[(1-e^{t})^{n-1}e^{t}\big]_{t=0}\nonumber\\
&=&\sum_{n=1}^{\infty}\frac{1}{n(n+1)}\frac{d^{k}}{dt^{k}}\big[(1-e^{t})^{n}\big]_{t=0}.\label{sec4-eq6}
\end{eqnarray}

But  $\frac{d^{k}}{dt^{k}}\big[(1-e^{t})^{n}\big]_{t=0}=0$ if
$n\ge k+1$. Therefore, the infinite sum in (\ref{sec4-eq6})
reduces to a finite sum

\begin{eqnarray}
I^{k-1}\big[\phi(-t)\big]_{t=0}&=&\sum_{n=1}^{k}\frac{1}{n(n+1)}\frac{d^{k}}{dt^{k}}\big[(1-e^{t})^{n}\big]_{t=0}\nonumber\\
&=&\sum_{n=1}^{k}\frac{1}{n(n+1)}\sum_{j=0}^{n}(-1)^{j}{n\choose
j}\frac{d^{k}}{dt^{k}}\big[e^{jt}\big]_{t=0}\nonumber\\
&=&(-1)^k\sum_{n=1}^{k}\frac{1}{n(n+1)}\sum_{j=0}^{n}(-1)^{j}{n\choose
j}j^{k}\nonumber\\
&=&
(-1)^k\sum_{n=1}^{k}\frac{1}{n(n+1)}\Delta_{n}(k).\label{sec4-eq7}
\end{eqnarray}

The last equation combined with (\ref{sec3-eq10}) gives the
explicit formula for $B_k$:

\begin{equation}\label{sec4-eq8}
B_{k}=(-1)^{k+1}\sum_{n=1}^{k}\frac{1}{n(n+1)}\Delta_{n}(k), k\ge
1.
\end{equation}

\section{Bernoulli Polynomials $B_k(1-x)$}\label{sec5}

The generalization of the series formula (\ref{sec4-eq1}) for the
Hurwitz zeta function $\zeta(s,x)$ defined for $0<x\le 1$ by

\begin{equation}\label{sec5-eq1}
\zeta(s,x)=\sum_{n=1}^{\infty} \frac{1}{(n-1+x)^{s}}
\end{equation}

is given in \cite{lazhar2} by the formula

\begin{equation}\label{sec5-eq2} (s-1)\zeta(s,x)
=\sum_{n=1}^{\infty}S_{n}(s,x)\bigg(\frac{1}{n+1}+
\frac{x-1}{n}\bigg),
\end{equation}

where  $S_{n}(s,x)$ is the generalization of the sums $S_n(s)$:

\begin{equation}\label{sec5-eq3}
S_{n}(s,x) = \sum_{k=0}^{n-1}(-1)^{k}{n-1\choose k}(k+x)^{-s}~{\rm
for}~ n\ge 2.
\end{equation}

There is also a corresponding integral given by
\begin{equation}\label{sec5-eq4}
(s-1)\zeta(s,x)
=\frac{1}{\Gamma(s)}\int_{0}^{\infty}\phi_x(t)t^{s-1}\,dt,
\end{equation}

where $\phi_x(t)$ is defined by

\begin{eqnarray}
\phi_x(t)&=&\frac{te^{-xt}}{(e^{t}-1)^2}-\frac{e^{-(x-1)t}}{e^{t}-1}+\frac{(x-1)te^{-(x-1)t}}{e^{t}-1}\nonumber \\
&=&\frac{d}{dt}\bigg(\frac{-t}{e^{t}-1}\bigg)e^{-(x-1)t}+\frac{(x-1)te^{-(x-1)t}}{e^{t}-1}\nonumber\\
&=&
\frac{d}{dt}\bigg(\frac{-te^{-(x-1)t}}{e^{t}-1}\bigg).\label{sec5-eq5}
\end{eqnarray}

The formula for Bernoulli polynomials is obtained by repeating
exactly the same steps of the previous section. We will repeat
these steps for the sake of clarity.

The Bernoulli polynomials are usually defined using the generating
function

\begin{equation}\label{sec5-eq6}
\frac{te^{xt}}{e^{t}-1}=\sum_{nk=0}^{\infty}\frac{B_k(x)}{k!}t^k,\quad
|t|<2\pi,
\end{equation}

or in terms of the function $\phi_x(t)$ as a generating function

\begin{equation}\label{sec5-eq7}
\phi_x(t)=\frac{d}{dt}\bigg(\frac{-te^{-(x-1)t}}{e^{t}-1}\bigg)=\sum_{k=1}^{\infty}\frac{-B_k(1-x)}{(k-1)!}t^{k-1}.
\end{equation}

Hence, by using the well-known identity $B_k(1-x)=(-1)^kB_k(x)$,
we finally obtain

\begin{equation}\label{sec5-eq8}
(-1)^kB_k(x)=B_k(1-x)=-\phi_x^{(k-1)}(0)=-I^{k-1}\big[\phi_x(-t)\big]_{t=0}.
\end{equation}

Since, for $0<t<\infty$,

\begin{equation}\label{sec5-eq9}
\phi_x(t)=\sum_{n=1}^{\infty}\bigg(\frac{1}{n+1}+\frac{x-1}{n}\bigg)(1-e^{-t})^{n-1}e^{-xt}
\end{equation}

uniformly, we may apply the fractional derivative operator
termwise to obtain

\begin{equation}\label{sec5-eq10}
I^{-s}\big[\phi_x(-t)\big]_{t=0}=\sum_{n=1}^{\infty}I^{-s}\bigg(\frac{1}{n+1}+\frac{x-1}{n}\bigg)\big[(1-e^{t})^{n-1}e^{xt}\big]_{t=0}.
\end{equation}

For $-s=k-1$,

\begin{eqnarray}
I^{k-1}\big[\phi_x(-t)\big]_{t=0}&=&\sum_{n=1}^{\infty}I^{k-1}\bigg(\frac{1}{n+1}+\frac{x-1}{n}\bigg)\big[(1-e^{t})^{n-1}e^{xt}\big]_{t=0}\nonumber\\
&=&\sum_{n=1}^{\infty}\bigg(\frac{1}{n+1}+\frac{x-1}{n}\bigg)\frac{d^{k-1}}{dt^{k-1}}\big[(1-e^{t})^{n-1}e^{xt}\big]_{t=0}\label{sec5-eq11}
\end{eqnarray}

But

\begin{eqnarray*}
&&\frac{d^{k-1}}{dt^{k-1}}\big[(1-e^{t})^{n-1}e^{-xt}\big]_{t=0}=0,~\rm{for}~ n\ge k+1 ~\rm{and}\\
&&(1-e^{t})^{n-1}e^{xt}=\sum_{j=0}^{n-1}(-1)^{j}{n-1\choose
k}e^{(j+x)t};
\end{eqnarray*}

therefore, the infinite sum in (\ref{sec5-eq6}) reduces to a
finite sum

\begin{eqnarray}
I^{k-1}\big[\phi_x(-t)\big]_{t=0}&=&\sum_{n=1}^{k}\bigg(\frac{1}{n+1}+\frac{x-1}{n}\bigg)\sum_{j=0}^{n-1}(-1)^{j}{n-1\choose
j}\frac{d^{k-1}}{dt^{k-1}}\big[e^{(j+x)t}\big]_{t=0}\nonumber\\
&=&(-1)^{k-1}\sum_{n=1}^{k}\bigg(\frac{1}{n+1}+\frac{x-1}{n}\bigg)\sum_{j=0}^{n-1}(-1)^{j}{n-1\choose
j}(j+x)^{k-1}\nonumber\\
&=&(-1)^{k-1}\sum_{n=1}^{k}\bigg(\frac{1}{n(n+1)}+\frac{x-1}{n^2}\bigg)\sum_{j=1}^{n}(-1)^{j-1}{n\choose
j}j(j+x-1)^{k-1}
\nonumber\\
&=&
(-1)^{k}\sum_{n=1}^{k}\bigg(\frac{1}{n(n+1)}+\frac{x-1}{n^2}\bigg)\Delta_{n,x}(k),\label{sec5-eq12}
\end{eqnarray}
where

\begin{equation}\label{sec5-eq13}
\Delta_{n,x}(k)=\sum_{j=1}^{n}(-1)^{j}{n\choose j}j(j+x-1)^{k-1}.
\end{equation}

The last equation combined with (\ref{sec5-eq8}) gives the
explicit formula for $B_k(1-x)$:

\begin{equation}\label{sec5-eq14}
B_{k}(1-x)=(-1)^{k+1}\sum_{n=1}^{k}\bigg(\frac{1}{n(n+1)}+\frac{x-1}{n^2}\bigg)\Delta_{n,x}(k),
k\ge 1.
\end{equation}

\section{Bernoulli Numbers and Polynomials  of Complex Index $s$}\label{sec6}

There are many generalizations of  integer-indexed Bernoulli
numbers and polynomials to  complex-indexed quantities. The reader
may consult for example \cite{butzer} or \cite{temme} and the
references therein. Here, we approach the generalization using
fractional derivatives.

The Bernoulli numbers and polynomials $B_s(x)$ for $s$ complex can
be defined using the fractional derivative operator of order $1-s$
(i.e. replace $k-1$ by $s-1$). When $\re(s)>0$ equation
(\ref{sec3-eq2}) applies. Formulas (\ref{sec4-eq3}) and
(\ref{sec4-eq4}) yield the following natural definition when
$\re(s)>0$:

\begin{eqnarray}
B_s=-I^{s-1}\big[\phi(-t)\big]_{t=0}&=&-\sum_{n=1}^{\infty}\frac{I^{s-1}\big[(1-e^{t})^{n-1}e^{t}\big]_{t=0}}{n+1}\nonumber\\
&=&-\sum_{n=1}^{\infty}\frac{S_{n}(1-s)}{n+1}\nonumber\\
&=&s\zeta(1-s).\label{sec6-eq1}
\end{eqnarray}

As for Bernoulli polynomials, their extension is obtained as
follows

\begin{eqnarray}
B_s(1-x)=-I^{s-1}\big[\phi_x(-t)\big]_{t=0}&=&-
\sum_{n=1}^{\infty}\bigg(\frac{1}{n+1}+\frac{x-1}{n}\bigg)I^{s-1}\big[(1-e^{t})^{n-1}e^{xt}\big]_{t=0}\nonumber\\
&=&-
\sum_{n=1}^{\infty}\bigg(\frac{1}{n+1}+\frac{x-1}{n}\bigg)I^{s-1}\big[\sum_{j=0}^{n-1}(-1)^{j}{n-1\choose
k}e^{(j+x)t}\big]_{t=0}\nonumber\\
&=&-\sum_{n=1}^{\infty}\frac{S_{n}(1-s,x)}{n+1}\nonumber\\
&=&s\zeta(1-s,x),\label{sec6-eq2}
\end{eqnarray}

Thus, from equation (\ref{sec6-eq2}), we see that the Bernoulli
polynomials extend naturally to the entire function
$s\zeta(1-s,x)$. This is an illustration that entire functions are
natural generalization of polynomials.

\section{The Second Explicit Formula for $B_k$}\label{sec7}

In this section we will prove formula (\ref{sec1-eq6}) using the
globally convergent series representation of the Riemann zeta
function discovered in \cite{hasse}. The series is given by

\begin{equation}\label{sec7-eq1}
s\zeta(s+1)=\sum_{n=0}^{\infty}\frac{S_{n+1}(s)}{n+1},
\end{equation}

$S_{n}(s)$ being defined in (\ref{sec4-eq2}).

Ii is easy to show that the sum $S_{n+1}(s)$ can be rewritten as
\begin{equation}\label{sec7-eq2}
S_{n+1}(s)=\frac{1}{\Gamma(s)}\int_{0}^{\infty}(1-e^{-t})^{n}e^{-t}t^{s-1}\,dt,
\end{equation}

and that for $\re(s)>0$,

\begin{equation}\label{sec7-eq3}
s\zeta(s+1)
=\frac{1}{\Gamma(s)}\int_{0}^{\infty}\psi(t)t^{s-1}\,dt,
\end{equation}

where the function $\psi(t)$ is given by
\begin{equation}\label{sec7-eq4}
\eta(t)=\frac{t}{e^{t}-1}=\sum_{n=0}^{\infty}\frac{(1-e^{-t})^{n}e^{-t}}{n+1}.
\end{equation}

Using the generating function $\psi(t)$, the Bernoulli numbers are
now given by

\begin{equation}\label{sec7-eq5}
B_k=\eta^{(k)}(0)=-I^{k}\big[\eta(-t)\big]_{t=0}.
\end{equation}

Again, we may apply the fractional derivative operator ($-s=k$)
termwise to obtain

\begin{eqnarray}
I^{k}\big[\eta(-t)\big]_{t=0}&=&\sum_{n=0}^{\infty}\frac{I^{k}\big[(1-e^{t})^{n}e^{t}\big]_{t=0}}{n+1}\nonumber\\
&=&\sum_{n=0}^{\infty}\frac{1}{n+1}\frac{d^{k}}{dt^{k}}\big[(1-e^{t})^{n}e^{t}\big]_{t=0}\nonumber\\
&=&\sum_{n=0}^{\infty}\frac{1}{(n+1)^2}\frac{d^{k+1}}{dt^{k+1}}\big[(1-e^{t})^{n+1}\big]_{t=0}.\label{sec7-eq6}
\end{eqnarray}

But  $\frac{d^{k+1}}{dt^{k+1}}\big[(1-e^{t})^{n+1}\big]_{t=0}=0$
if $n+1\ge k+2$. Therefore, the infinite sum in (\ref{sec7-eq6})
reduces to a finite sum

\begin{eqnarray}
I^{k}\big[\eta(-t)\big]_{t=0}&=&\sum_{n=0}^{k}\frac{1}{(n+1)^2}\frac{d^{k+1}}{dt^{k+1}}\big[(1-e^{t})^{n+1}\big]_{t=0}\nonumber\\
&=&\sum_{n=0}^{k}\frac{1}{(n+1)^2}\sum_{j=0}^{n+1}(-1)^{j}{n+1\choose
j}\frac{d^{k+1}}{dt^{k+1}}\big[e^{jt}\big]_{t=0}\nonumber\\
&=&(-1)^{k+1}\sum_{n=0}^{k}\frac{1}{(n+1)^2}\sum_{j=0}^{n+1}(-1)^{j}{n+1\choose
j}j^{k+1}\nonumber\\
&=&
(-1)^{k+1}\sum_{n=0}^{k}\frac{1}{(n+1)^2}\Delta_{n+1}(k+1).\label{sec7-eq7}
\end{eqnarray}

With an appropriate change of variable in the summation index, the
last equation combined with (\ref{sec7-eq5}) gives the explicit
formula for $B_k$:

\begin{equation}\label{sec7-eq8}
B_{k}=(-1)^{k+1}\sum_{n=1}^{k+1}\frac{1}{n^2}\Delta_{n}(k+1), k\ge
0.
\end{equation}

The extension of the last explicit formula to Bernoulli
polynomials and Bernoulli numbers of fractional index is
straightforward.





\bibliographystyle{model1-num-names}

\end{document}